# Construction of Additive Semi-Implicit Runge-Kutta methods with low-storage requirements*

I. Higueras and T. Roldán

Dpto. Ingeniería Matemática e Informática
Universidad Pública de Navarra
higueras@unavarra.es
teo@unavarra.es

## Abstract

Space discretization of some time-dependent partial differential equations gives rise to systems of ordinary differential equations in additive form whose terms have different stiffness properties. In these cases, implicit methods should be used to integrate the stiff terms while efficient explicit methods can be used for the non-stiff part of the problem. However, for systems with a large number of equations, memory storage requirement is also an important issue. When the high dimension of the problem compromises the computer memory capacity, it is important to incorporate low memory usage to some other properties of the scheme.

In this paper we consider Additive Semi-Implicit Runge-Kutta (ASIRK) methods, a class of implicit-explicit Runge-Kutta methods for additive differential systems. We construct two second order 3-stage ASIRK schemes with low-storage requirements. Having in mind problems with stiffness parameters, besides accuracy and stability properties, we also impose stiff accuracy conditions. The numerical experiments done show the advantages of the new methods.

## 1 Introduction

Space discretization of some time-dependent partial differential equations (PDEs) gives rise to systems of ordinary differential equations in additive form

$$y' = f(y) + g(y), \qquad y(t_0) = y_0, \tag{1}$$

where $f, g : \mathbb{R}^N \to \mathbb{R}^N$ are sufficiently smooth functions with different stiffness properties. These systems, often with a large number of equations, arise, e.g., from semidiscretisations of convection-diffusion problems and hyperbolic systems with relaxation [4, 5, 16, 22, 23].

When a differential system (1) involves terms with different stiffness properties, a natural approach to obtain numerical approximations is the use of implicit-explicit (IMEX) time-discretizations; in this procedure, an explicit method is used for the non-stiff part $f$, and an implicit one is used for the stiff part $g$. IMEX Runge-Kutta methods have been deeply studied in the literature (see, e.g., [1, 3, 4, 16, 22, 28]).

Sometimes, the differential system (1) is of the form

$$\begin{aligned} u' &= f_1(u, v), \\ v' &= f_2(u, v) + \frac{1}{\varepsilon} g_2(u, v), \end{aligned} \tag{2}$$

---

*Supported by Ministerio de Economía y Competitividad, project MTM2011-23203.

where $\varepsilon$ is the stiffness parameter; observe that (2) is a particular case of (1) with $y = (u,v)$, $f = (f_1, f_2)$, $g = (0, g_2/\varepsilon)$. These systems have been considered, e.g., in [1, 2, 3, 4, 22, 23], where robust IMEX Runge-Kutta methods have been analyzed. In particular, in [3, 4], uniform convergence in the stiffness parameter $\varepsilon$ has been studied.

For systems with a large number of equations, memory storage requirement is an important issue. When the high dimension of the problem compromises the computer memory capacity, it is important to incorporate low memory usage to some other properties of the scheme. These ideas have been developed in [7, 6, 8, 10, 13, 17, 18, 19, 20, 26, 27], where different low-storage Runge-Kutta methods have been constructed. In particular, in [7, 17], explicit Runge-Kutta methods are implemented by using the van der Houwen format [26], while diagonally implicit Runge-Kutta (DIRK) schemes have been explored in [20].

The aim of this paper is to construct robust IMEX Runge-Kutta methods with low-storage requirements. We consider Additive Semi-Implicit Runge-Kutta methods (ASIRK), a special class of IMEX Runge-Kutta methods considered, e.g., in [28]. Having in mind problems of the form (2), on the new schemes we impose, not only accuracy and stability properties but also uniform accuracy in the stiffness parameter $\varepsilon$.

The rest of the paper is organized as follows. In section 2 we introduce ASIRK methods, and we study how to incorporate to them accuracy, stability and stiff accuracy properties. Section 3 is devoted to the study of low-storage implementations of ASIRK schemes; from the coefficient pattern of low-storage explicit Runge-Kutta and DIRK methods, we will get the coefficient form of low-storage ASIRK schemes. With the results of the previous sections, two ASIRK methods are constructed in section 4. Some numerical experiments done in section 5 show the advantages of the new schemes. The paper ends with some concluding remarks.

## 2 A class of implicit-explicit Runge-Kutta methods

In this section we consider a special class of IMEX Runge-Kutta schemes. In the first part of the section, we review some known results on additive Runge-Kutta schemes. Next, we formulate ASIRK schemes as additive Runge-Kutta methods and we obtain the set of order conditions to obtain third order of accuracy. The section ends with a review of linear stability for additive Runge-Kutta methods, and a stiff accuracy analysis for a two-dimensional model of the form (2).

**Additive Runge-Kutta methods**

We consider $s$-stage additive Runge-Kutta methods with coefficients $(\tilde{A}, \tilde{b})$, $(A, b)$. If $y_n$ denotes the numerical approximation of the solution of (1) at time $t_n = t_0 + nh$, then the numerical approximation $y_{n+1}$ to $y(t_{n+1})$ with an $s$-stage additive Runge-Kutta method is given by

$$y_{n+1} = y_n + h \sum_{i=1}^{s} \tilde{b}_i f(Y_{i,n+1}) + h \sum_{i=1}^{s} b_i g(Y_{i,n+1}), \tag{3}$$

where the internal stages $Y_{i,n+1}$ are obtained from

$$Y_{i,n+1} = y_n + h \sum_{j=1}^{s} \tilde{a}_{ij} f(Y_{j,n+1}) + h \sum_{j=1}^{s} a_{ij} g(Y_{j,n+1}), \quad i = 1, \ldots, s.$$

The coefficients of the method are usually represented in a double Butcher tableau

$$\begin{array}{c|c} \tilde{c} & \tilde{A} \\ \hline & \tilde{b}^t \end{array}, \quad \begin{array}{c|c} c & A \\ \hline & b^t \end{array}, \tag{4}$$

where $\tilde{c}_i = \sum_{j=1}^{s} \tilde{a}_{ij}$ and $c_i = \sum_{j=1}^{s} a_{ij}$, $i = 1, \ldots, s$.

The equation (1) can be displayed as a partitioned system

$$u' = \hat{f}(u,v),$$
$$v' = \hat{g}(u,v),$$

with $y = u + v$, $\hat{f}(u,v) = f(u+v)$ and $\hat{g}(u,v) = g(u+v)$. Then the additive Runge-Kutta method can be seen as a particular class of partitioned Runge-Kutta scheme, and the order conditions can be obtained from the theory of order conditions for partitioned methods [12]. In particular the 20 order conditions (see [16] for details) up to order $p = 3$ can be simplified to the following ones if the abscissa vectors $\tilde{c}$ and $c$ coincide.

$$
\begin{array}{lll}
\text{First order:} & \tilde{b}^t e = 1 & b^t e = 1 \\
\text{Second order:} & \tilde{b}^t c = 1/2 & b^t c = 1/2 \\
\text{Third order:} & \tilde{b}^t \tilde{A} c = 1/6 & b^t A c = 1/6 \\
& \tilde{b}^t c^2 = 1/3 & b^t c^2 = 1/3 \\
& \tilde{b}^t A c = 1/6 & b^t \tilde{A} c = 1/6
\end{array}
\tag{5}
$$

where, as usual, $e = (1, \ldots, 1)^t \in \mathbb{R}^s$.

When system (1) contains terms with different stiffness properties, different kinds of methods can be used for each term. For example, if $f$ is the non-stiff part and $g$ is the stiff one, an explicit method can be used for $f$ and an implicit one for $g$. In this way, we obtain IMEX Runge-Kutta methods. Furthermore, aiming at diminishing the computational cost, it is usual to consider DIRK schemes for the implicit part. In this case, the internal stages $Y_{i,n+1}$, $i = 1, \ldots s$, are given by

$$Y_{i,n+1} = y_n + h \sum_{j=1}^{i-1} \tilde{a}_{ij} f(Y_{j,n+1}) + h \sum_{j=1}^{i} a_{ij} g(Y_{j,n+1}), \quad i = 1, \ldots, s. \tag{6}$$

This kind of IMEX Runge-Kutta schemes has been used in [1, 2, 3, 4, 22, 23] in different applications.

**ASIRK-$s$A methods as additive Runge-Kutta methods**

In [14, 16, 28], looking for computational efficiency, different additive semi-implicit Runge-Kutta methods are considered. In this paper we are interested in a class of additive methods studied in [28], namely the ones named ASIRK-$s$A. The numerical solution of (1) with an ASIRK-$s$A method is given by

$$y_{n+1} = y_n + \sum_{i=1}^{s} \omega_i K_{i,n+1}, \tag{7}$$

where the internal derivatives $K_{i,n+1}$ are given by

$$K_{i,n+1} = h\left(f\left(y_n + \sum_{j=1}^{i-1} b_{ij} K_{j,n+1}\right) + g\left(y_n + \sum_{j=1}^{i-1} c_{ij} K_{j,n+1} + c_{ii} K_{i,n+1}\right)\right), \quad i = 1, \ldots, s, \tag{8}$$

and $b_{ij}$, $c_{ij}$, $\omega_j$ are the coefficients of the method. Observe that the resulting scheme is explicit in $f$ and diagonally implicit in $g$. Some other schemes, like ASIRK-$s$B and ASIRK-$s$C methods, are also considered in [28]. These schemes are semi-implicit extension of Rosenbrock Runge–Kutta methods

$$(I_N - hc_{ii}J(y_n + \sum_{j=1}^{i-1} d_{ij} K_{j,n+1}))K_{i,n+1} = h\left(f\left(y_n + \sum_{j=1}^{i-1} b_{ij} K_{j,n+1}\right) + g\left(y_n + \sum_{j=1}^{i-1} c_{ij} K_{j,n+1}\right)\right),$$
$$i = 1, \ldots, s,$$

where $J = \partial g / \partial y$ and $d_{ij}$ are additional parameters. Although these schemes are computationally more efficient, stability reasons make ASIRK-$s$A methods more interesting.

On the following, we denote the matrices and vector containing the coefficients of the ASIRK-$s$A method by $\mathscr{B} = (b_{ij})$, $\mathscr{C} = (c_{ij})$, and $\omega$, respectively. Below we display these matrices for $s = 3$,

$$\mathscr{B} = \begin{pmatrix} 0 & 0 & 0 \\ b_{21} & 0 & 0 \\ b_{31} & b_{32} & 0 \end{pmatrix}, \quad \mathscr{C} = \begin{pmatrix} c_{11} & 0 & 0 \\ c_{21} & c_{22} & 0 \\ c_{31} & c_{32} & c_{33} \end{pmatrix}, \quad \omega = \begin{pmatrix} \omega_1 \\ \omega_2 \\ \omega_3 \end{pmatrix}. \tag{9}$$

Even though method (7)-(8) can be considered as an additive Runge-Kutta method, it is important to notice that the coefficients in (9) are not the standard coefficients (4) of an additive Runge-Kutta scheme (3)-(6). However, as we show below, equations (7)-(8) can be rewritten as an additive Runge-Kutta method (3)-(6). With this approach, some issues for ASIRK-$s$A methods, e.g., the set of order conditions, can be obtained in a simple way from the theory of additive Runge-Kutta methods.

Equations (7)-(8) are given in terms of $K_{n+1} = (K_{1,n+1}^t, \ldots, K_{s,n+1}^t)^t$, the internal derivative vector. By using the Kronecker product $\otimes$ and denoting $K_{n+1} = h\left(F(Y_{n+1}) + G(\hat{Y}_{n+1})\right)$, where $F(Y_{n+1}) = (f(Y_{1,n+1})^t, \ldots, f(Y_{s,n+1})^t)^t$ and $G(\hat{Y}_{n+1}) = (g(\hat{Y}_{1,n+1})^t, \ldots, g(\hat{Y}_{s,n+1})^t)^t$, it is possible to rewrite (7)-(8) as

$$y_{n+1} = y_n + h\left(\omega^t \otimes I_N\right)\left(F(Y_{n+1}) + G(\hat{Y}_{n+1})\right), \tag{10}$$

where the internal stages $Y_{n+1} = (Y_{1,n+1}^t, \ldots, Y_{s,n+1}^t)^t$ and $\hat{Y}_{n+1} = (\hat{Y}_{1,n+1}^t, \ldots, \hat{Y}_{s,n+1}^t)^t$ are obtained from

$$Y_{n+1} = e \otimes y_n + (\mathscr{B} \otimes I_N)\left(h F(Y_{n+1}) + h G(\hat{Y}_{n+1})\right), \tag{11}$$

$$\hat{Y}_{n+1} = e \otimes y_n + (\mathscr{C} \otimes I_N)\left(h F(Y_{n+1}) + h G(\hat{Y}_{n+1})\right). \tag{12}$$

This formulation corresponds to a $2s$-stage additive Runge-Kutta method, where the $2s$ internal stages are $Y_{n+1}$ and $\hat{Y}_{n+1}$. Although the natural computation of the internal stages (11)-(12) is $(Y_{1,n+1}^t, \hat{Y}_{1,n+1}^t, \ldots, Y_{s,n+1}^t, \hat{Y}_{s,n+1}^t)^t$, we have preferred to write them as $(Y_{1,n+1}^t, \ldots, Y_{s,n+1}^t, \hat{Y}_{1,n+1}^t, \ldots, \hat{Y}_{s,n+1}^t)^t$ because the corresponding equivalent additive Runge-Kutta method (4) can be written with the following block structure

$$\begin{array}{c|cc|cc} \mathscr{B}e & \mathscr{B} & 0 & 0 & \mathscr{B} \\ \mathscr{C}e & \mathscr{C} & 0 & 0 & \mathscr{C} \\ \hline & \omega^t & 0 & 0 & \omega^t \end{array} \tag{13}$$

and it is easy to obtain the set of order conditions from (13). Observe that (13) is an IMEX Runge-Kutta method (4) with $\tilde{c} = c$. Now, from the set of order conditions for IMEX methods (5) we can derive the set of order conditions for ASIRK-$s$A scheme (7)-(8) in terms of matrices $\mathscr{B}$, $\mathscr{C}$, and vector $\omega$. Below we display these equations, particularized for $s = 3$, up to order $p = 3$.

First order:

$$\omega^t e = 1 \qquad\qquad \omega_1 + \omega_2 + \omega_3 = 1 \tag{14}$$

Second order:

$$\omega^t \mathscr{B} e = 1/2 \qquad\qquad \omega_2 b_{21} + \omega_3 (b_{31} + b_{32}) = 1/2 \tag{15}$$

$$\omega^t \mathscr{C} e = 1/2 \qquad\qquad \omega_1 c_{11} + \omega_2 (c_{22} + c_{21}) + \omega_3 (c_{33} + c_{31} + c_{32}) = 1/2 \tag{16}$$

Third order:

$$\omega^t \mathscr{B}^2 e = 1/6 \qquad\qquad \omega_3 b_{21} b_{32} = 1/6 \tag{17}$$

$$\omega^t (\mathscr{B}e)^2 = 1/3 \qquad\qquad \omega_2 b_{21}^2 + \omega_3 (b_{31} + b_{32})^2 = 1/3 \tag{18}$$

$$\omega^t \mathscr{C}^2 e = 1/6 \qquad\qquad \omega_1 c_{11}^2 + \omega_2 (c_{21} c_{11} + c_{22} c_{21} + c_{22}^2) + \tag{19}$$
$$\omega_3 (c_{31} c_{11} + c_{33}^2 + c_{33} c_{31} + c_{22} c_{32} + c_{33} c_{32} + c_{21} c_{32}) = 1/6$$

$$\omega^t (\mathscr{C}e)^2 = 1/3 \qquad\qquad \omega_1 c_{11}^2 + \omega_2 (c_{22} + c_{21})^2 + \omega_3 (c_{33} + c_{31} + c_{32})^2 = 1/3 \tag{20}$$

$$\omega^t \mathscr{B}\mathscr{C} e = 1/6 \qquad\qquad c_{11} (\omega_2 b_{21} + \omega_3 b_{31}) + \omega_3 b_{32} (c_{22} + c_{21}) = 1/6 \tag{21}$$

$$\omega^t \mathscr{C}\mathscr{B} e = 1/6 \qquad\qquad \omega_3 c_{33} (b_{31} + b_{32}) + b_{21} (\omega_2 c_{22} + \omega_3 c_{32}) = 1/6 \tag{22}$$

Observe that (17)-(18) and (19)-(20) are order-3 conditions for the explicit and the implicit method, respectively. However, (21)-(22) are coupling conditions, including coefficients $b_{ij}$ and $c_{ij}$ from both, the explicit and the implicit scheme.

For some additive problems (1) it may occur that $f_y g_y = g_y f_y$. In such case, a simple analysis shows that equations (21)-(22) are reduced to the sum of both of them,

$$\omega^t \mathscr{B}\mathscr{C} e + \omega^t \mathscr{C}\mathscr{B} e = 1/3. \qquad (23)$$

**Linear stability for ASIRK-$s$A methods**

In this section we study linear stability properties of ASIRK-$s$A methods (7)-(8). To derive them, we apply some known results on linear stability for additive Runge-Kutta methods [16, 22, 28]; our approach is similar to the one made in [22].

We consider the simplified linear model equation

$$y' = \xi_1 y + \xi_2 y, \qquad y(0) = 1, \qquad (24)$$

where $\xi_1, \xi_2 \in \mathbb{C}$, with $\operatorname{Re}(\xi_2) < 0$. If we solve the model equation (24) with method (7)-(8), the numerical solution $y_1$ after one step $h$ is given by

$$y_1 = R(z_1, z_2) = 1 + (z_1 + z_2)\omega^t (I_s - z_1 \mathscr{B} - z_2 \mathscr{C})^{-1} e \qquad (25)$$

where $z_1 = \xi_1 h$, $z_2 = \xi_2 h$. The expression $R(z_1, z_2)$ is the function of absolute stability, and

$$S_A = \{(z_1, z_2) \in \mathbb{C}^2 : |R(z_1, z_2)| \leq 1\}$$

is the region of absolute stability. We are interested in computing the largest set

$$S_1 = \{z_1 \in \mathbb{C} : \sup_{z_2 \in \mathbb{C}^-} |R(z_1, z_2)| \leq 1\}. \qquad (26)$$

For more details, see [22]. Besides, in order to get a correct asymptotic decay for the stiff terms, the $L$-stability condition for the implicit step, that is,

$$\lim_{z_2 \to -\infty} R(0, z_2) = 1 - \omega^t \mathscr{C}^{-1} e = 0, \qquad (27)$$

must be imposed.

In section 4, in the process of ASIRK-$s$A method construction, we will compute the set $S_1$ and we will impose condition (27) to the implicit scheme. In this way, if a method is $A$-stable, it will also be $L$-stable.

**Stiff accuracy analysis**

In order to increase the accuracy of the schemes constructed in this paper, we will impose stiff accuracy, that is, accuracy of the numerical solution for problems of the form (2) for small values of $\varepsilon$. This issue has been studied in [4] and [22] for IMEX Runge-Kutta methods. In [4], the analysis done leads to the additional order conditions

$$b^t A^{-1} \tilde{c} = 1, \quad b^t A^{-1} \tilde{c}^2 = 1, \quad b^t A^{-1} \tilde{A} \tilde{c} = 1/2, \qquad (28)$$

that ensure accuracy also in the stiff regime when the parameter $\varepsilon$ tends to zero.

In our case, when an ASIRK-$s$A method is written as an IMEX Runge-Kutta scheme, matrix $A$ in (13) is not invertible. In order to impose conditions (28), we introduce a $\rho$ parameter in the formulation (13)

$$\begin{array}{c|cc|cc}
\mathscr{B}e & \mathscr{B} & 0 & \rho I_s & \mathscr{B} \\
\mathscr{C}e & \mathscr{C} & 0 & 0 & \mathscr{C} \\
\hline
 & \omega^t & 0 & 0 & \omega^t
\end{array}$$

and apply conditions (28). Then, in the limit case $\rho \to 0$, we get

$$\omega^t e = 1, \quad \omega^t \mathscr{C}^{-1} (\mathscr{C} e)^2 = 1, \quad \omega^t \mathscr{B} e = 1/2. \quad (29)$$

Observe that the first and the third equations are order conditions (14) and (15), respectively. Furthermore, if the last row of $\mathscr{C}$ is $\omega$, as $\omega^t \mathscr{C}^{-1} = (0, \ldots, 0, 1)^t$, the condition in the middle reduces to $(\omega^t e)^2 = 1$, and thus it is fulfilled for first order schemes.

Following [22], in our study we will also consider the linear system

$$\begin{cases} u' = \delta_1 u + \sigma_1 v, \\ v' = \delta_2 u + \sigma_2 v + \dfrac{1}{\varepsilon}(c u - v), \end{cases} \quad (30)$$

as a particular case of (2). This problem has been analyzed in [22], where the authors attempt to obtain robust IMEX methods for hyperbolic systems with relaxation. If we consider consistent initial values, $v_0 = c u_0$, then, after one time step, the exact solution of (30) satisfies

$$u(\varepsilon, h) = u_0 \left(1 + \hat{a} h + \tfrac{1}{2} \hat{a}^2 h^2 + \sigma_1 \hat{b} h \varepsilon\right) + \mathcal{O}(h^3, h \varepsilon^2),$$

$$v(\varepsilon, h) = u_0 \left((1 + \hat{a} h + \tfrac{1}{2} \hat{a}^2 h^2) c + \hat{b} \varepsilon + (\hat{a} + \sigma_1 c) \hat{b} h \varepsilon\right) + \mathcal{O}(h^3, h \varepsilon^2),$$

where, as in [22], we have denoted $\hat{a} = \delta_1 + \sigma_1 c$, and $\hat{b} = \delta_2 + (\sigma_2 - \delta_1) c - \sigma_1 c^2$.

In section 5 we will obtain the numerical solution $(u_1, v_1)$ of system (30) with different ASIRK-$s$A additive schemes (7)-(8). Then, from the difference between the numerical and the exact solution

$$u_1(\varepsilon, h) - u(\varepsilon, h), \qquad v_1(\varepsilon, h) - v(\varepsilon, h), \quad (31)$$

we will get conditions on the coefficients of the method.

## 3 Low-storage ASIRK-$s$A methods

As it has been pointed out in the introduction, for systems with a large number of equations, the efficiency of the numerical simulation increases if low-storage schemes are used. In order to get these low-storage implementations, the coefficient matrices of the methods must have a special pattern. In this section, we will study this issue for ASIRK-$s$A schemes. We will count the number of registers of memory required; each memory register is a vector of dimension $N$, where $N$ is the number of equations of the problem.

In general, for an ASIRK-$s$A method, $s+1$ memory registers are needed. For example, for $s = 3$, a basic implementation uses four memory registers: one for each stage derivative $K_{i,n+1}$, $i = 1, 2, 3$, and another one for the numerical solution $y_{n+1}$. However, as we will show later on, if we impose some conditions on the coefficients, implementations with less memory usage are possible.

As we have pointed out above, ASIRK-$s$A methods can be formulated as IMEX Runge-Kutta methods. However, as far as we know, low-storage aspects for this class of methods have not been deeply studied in the literature. Some comments about this topic are done in [13], where several strongly stable IMEX Runge-Kutta methods are implemented for the ANTARES code. In particular, though the authors do not claim that this is the best way to proceed, six vectors in memory are used for a 2-stage scheme, and seven vectors in memory are used for a 3-stage one.

In order to get low-storage implementations for ASIRK-$s$A methods, we first study low-storage explicit and DIRK Runge-Kutta schemes separately. In both cases, the internal stages are computed sequentially. In particular, we focus our attention on methods that require just two memory registers.

In the rest of the section, in order to simplify the notation, vectors $K_{i,n+1}$ and $Y_{i,n+1}$ will be denoted by $K_i$ and $Y_i$, respectively.

## 3.1 Explicit Runge-Kutta methods implemented in two memory registers

In this section we consider explicit Runge-Kutta methods; the equations for these methods can be directly obtained from the formulas (3) and (6) by just considering $g = 0$.

There are many papers on low-storage strategies for explicit Runge-Kutta schemes in the literature (see, e.g., [7, 6, 8, 10, 17, 18, 19, 26, 27]). In [7, 17], explicit Runge-Kutta methods are implemented by using the van der Houwen [26] procedure that uses 2 memory registers. An $s$-stage explicit Runge-Kutta method $(\tilde{A}, \tilde{b})$ implemented in this way is defined by the matrices

$$\tilde{A} = \begin{pmatrix} 0 & & & & \\ b_1+\gamma_1 & 0 & & & \\ b_1 & b_2+\gamma_2 & 0 & & \\ \vdots & \vdots & \ddots & \ddots & \\ b_1 & b_2 & \ldots & b_{s-1}+\gamma_{s-1} & 0 \end{pmatrix}, \quad \tilde{b} = \begin{pmatrix} b_1 \\ b_2 \\ \vdots \\ b_s \end{pmatrix}. \tag{32}$$

This way to proceed implies severe restrictions on the number of free coefficients of the method. For an $s$-stage method, there are only $2s - 1$ free parameters.

For example, the equations of a low-storage 3-stage explicit Runge-Kutta scheme (32) are

$$\begin{aligned} K_1 &= hf(y_n), \\ K_2 &= hf(y_n + b_1 K_1 + \gamma_1 K_1), \\ K_3 &= hf(y_n + b_1 K_1 + b_2 K_2 + \gamma_2 K_2), \\ y_{n+1} &= y_n + b_1 K_1 + b_2 K_2 + b_3 K_3, \end{aligned} \tag{33}$$

where the stage derivatives $K_i$ are computed sequentially. For methods of the form (32), for each stage, a memory register (Register 1) is used for the storage of the auxiliary value $y_n + \sum_{j=1}^{i-1} b_j K_j$,

$$\text{Register 1:} \quad Y_i = Y_{i-1} + b_{i-1} K_{i-1}, \quad i = 1, \ldots, s+1, \tag{34}$$

where we consider $b_0 = 0$ and $Y_0 = y_n$. After the last stage, we obtain the numerical solution $y_{n+1}$ in the Register 1 as $Y_{s+1} = Y_s + b_s K_s$. The second memory register (Register 2) is used for the storage of each stage derivative $K_i$

$$\text{Register 2:} \quad K_i = hf(Y_i + \gamma_{i-1} K_{i-1}), \quad i = 1, \ldots, s, \tag{35}$$

where we set $\gamma_0 = 0$.

## 3.2 DIRK methods implemented in two memory registers

Next, we study low-storage DIRK methods. The formulation of these methods can be obtained from the equations (3) and (6) by just assuming $f = 0$.

As far as we know, low-storage issues for DIRK schemes have not been seriously studied in the literature. In [20], singly diagonally implicit Runge-Kutta (SDIRK) methods with only two (three) vectors in memory for a second (third) order scheme are considered. As we will show in this section, imposing some conditions on the coefficients, DIRK schemes can be implemented with two memory registers for any number of stages.

**Proposition 1.** *Let $(A, b)$ be the coefficients of an $s$-stage DIRK method. Assume that the elements of $A$ and $b$, denoted by $a_{ij}$ and $b_j$, respectively, satisfy the following relationships*

$$\begin{aligned} a_{i+1,i} &= a_{i+2,i} = \cdots = a_{s,i} = b_i, \quad i = 1, \ldots, s-1, \\ \frac{a_{i+1,i}}{a_{i,i}} &= \alpha_i, \quad i = 1, \ldots, s-1, \\ \frac{b_s}{a_{ss}} &= \alpha_s, \end{aligned} \tag{36}$$

for some $\alpha_i \in \mathbb{R}$, $i = 1, \ldots, s$. Then, the DIRK method can be expressed as

$$Y_i = W_{i-1} + h a_{ii} g(Y_i), \qquad W_i = \alpha_i Y_i + \beta_i W_{i-1}, \qquad i = 1, \ldots, s, \tag{37}$$

where $W_0 = y_n$, $W_s = y_{n+1}$, and $\beta_i$ are defined by

$$\beta_i = 1 - \alpha_i, \qquad i = 1, \ldots, s. \tag{38}$$

*Conversely, if an s-stage DIRK method can be expressed as (37), then its coefficients must satisfy equations (36) and (38).*

*Proof.* The proof is straightforward. Following an induction process, we first check that, defining $W_0 = y_n$, from (37), we get

$$W_i = y_n + h \sum_{j=1}^{i} a_{i+1,j} g(Y_j), \qquad i = 1, \ldots, s-1.$$

In this way, the expression $W_{i-1} + h a_{ii} g(Y_i)$ is the $i$-th internal stage $Y_i$. Secondly, we check that $W_s = y_n + h \sum_{i=1}^{s} b_j g(Y_j)$.

To prove the second part, we simply have to check that, if $Y_i$ in (37) is the standard $i$-th internal stage of a DIRK method, $i = 1, \ldots, s$, then equations (36) and (38) must hold. □ □

DIRK schemes satisfying the conditions in Proposition 1 can be implemented using 2 memory registers in the following way. For each stage, a memory register (Register 1) is used for the storage of the internal stage $Y_i$,

$$\text{Register 1:} \quad Y_i = W_{i-1} + h a_{ii} g(Y_i), \quad i = 1, \ldots, s, \tag{39}$$

where $W_i$ is an auxiliary vector, and $W_0 = y_n$. The other memory register (Register 2) is used for the storage of this auxiliary vector

$$\text{Register 2:} \quad W_i = \alpha_i Y_i + \beta_i W_{i-1}, \quad i = 1, \ldots, s. \tag{40}$$

After the last stage, we obtain the numerical solution $y_{n+1}$ in this register as $W_s$.

Implementation of a DIRK method with two memory registers restricts the number of free coefficients of the method. For $s = 3$, the Butcher tableau is the following one

$$\begin{array}{c|ccc} c_1 & b_1/\alpha_1 & 0 & 0 \\ c_2 & b_1 & b_2/\alpha_2 & 0 \\ c_3 & b_1 & b_2 & b_3/\alpha_3 \\ \hline & b_1 & b_2 & b_3 \end{array} \tag{41}$$

**Proposition 1.** *In [9], in the context of Strong Stability Preserving methods, the optimal second order 3-stage SDIRK method is of the form (41) with coefficients $b_1 = b_2 = b_3 = 1/3$ and $\alpha_1 = \alpha_2 = \alpha_3 = 2$.*

**Proposition 2.** *If the differential problem is nonlinear, DIRK methods require the resolution of nonlinear system (39). If fixed point iteration is used to solve this system, no additional memory registers are needed. However, with other iterative schemes, e.g., Newton like methods, additional memory registers must be used. The study of nonlinear solvers and the number of memory registers required is beyond the aims of this paper and it is not considered here.*

**Proposition 3.** *If $\alpha_s = 1$, we obtain a stiffly accurate DIRK method, that is, the last row of A is equal to vector b. In this case, A-stable DIRK methods are also L-stable schemes.*

### 3.3 ASIRK-*s*A methods implemented in three memory registers

We are in position to propose a kind of ASIRK-*s*A methods that can be implemented in an efficient way by using three memory registers for any number of stages. The idea is to consider the form of scheme (33) for the explicit part of the ASIRK-*s*A method, and the form of scheme (41) for the implicit part. In this way, if function $g$ in (7)-(8) is equal to zero, we recover the van der Houwen implementation (34)-(35) for an explicit method; in a similar way, if function $f$ in (7)-(8) is equal to zero, we recover the implementation (39)-(40) for the implicit scheme.

On the following, for simplicity, we will write $\lambda_i$ in the diagonal entries instead of $\omega_i/\alpha_i$, $i=1,\ldots,s$. The ASIRK-*s*A methods we propose are of the form

$$\mathscr{B} = \begin{pmatrix} 0 & & & & \\ \omega_1+\gamma_1 & 0 & & & \\ \omega_1 & \omega_2+\gamma_2 & 0 & & \\ \vdots & \vdots & \ddots & \ddots & \\ \omega_1 & \omega_2 & \ldots & \omega_{s-1}+\gamma_{s-1} & 0 \end{pmatrix}, \quad \mathscr{C} = \begin{pmatrix} \lambda_1 & & & & \\ \omega_1 & \lambda_2 & & & \\ \omega_1 & \omega_2 & \lambda_3 & & \\ \vdots & \vdots & \ddots & \ddots & \\ \omega_1 & \omega_2 & \ldots & \omega_{s-1} & \lambda_s \end{pmatrix}, \quad \omega = \begin{pmatrix} \omega_1 \\ \omega_2 \\ \vdots \\ \omega_s \end{pmatrix}. \quad (42)$$

For example, an ASIRK-*s*A method of the form (42) for $s=3$ is

$$K_1 = hf(y_n) + hg(y_n + \lambda_1 K_1),$$
$$K_2 = hf(y_n + \omega_1 K_1 + \gamma_1 K_1) + hg(y_n + \omega_1 K_1 + \lambda_2 K_2),$$
$$K_3 = hf(y_n + \omega_1 K_1 + \omega_2 K_2 + \gamma_2 K_2) + hg(y_n + \omega_1 K_1 + \omega_2 K_2 + \lambda_3 K_3),$$
$$y_{n+1} = y_n + \omega_1 K_1 + \omega_2 K_2 + \omega_3 K_3.$$

It is possible to implement an *s*-stage ASIRK method of the form (42) by using just 3 memory vectors; we simply have to combine the process followed for explicit and DIRK methods. For each stage, a memory register (Register 1) is used for the storage of $y_n + \sum_{j=1}^{i-1}\omega_j K_j$,

$$\text{Register 1:} \quad Y_i = Y_{i-1} + \omega_{i-1}K_{i-1}, \quad i=1,\ldots,s+1,$$

where we consider $b_0=0$ and $Y_0 = y_n$. After the last stage, we obtain the numerical solution $y_{n+1}$ in the Register 1 as $Y_{s+1} = Y_s + \omega_s K_s$. The second memory register (Register 2) is used for the storage of the evaluation of the function $f$,

$$\text{Register 2:} \quad L_i = hf(Y_i + \gamma_{i-1}K_{i-1}), \quad i=1,\ldots,s,$$

where $\gamma_0 = 0$. Finally, the third memory register (Register 3) is used for the internal derivative $K_i$

$$\text{Register 3:} \quad K_i = L_i + hg(Y_i + \lambda_i K_i), \quad i=1,\ldots,s. \quad (43)$$

**Proposition 4.** *Observe that, if $\lambda_s = \omega_s$, the implicit scheme is stiffly accurate. Besides, in this case, the implicit part of the method satisfies condition* (27).

## 4 New ASIRK-3A methods implemented in three registers of memory

Once the structure (42) of low-storage ASIRK-*s*A methods is determined, we begin the construction of schemes imposing accuracy, stiff accuracy and stability properties. We will consider the cases $s=2,3$.

*Case $s=2$.* First, we set $\lambda_2 = \omega_2$ (see Remark (4)). Besides, we impose condition (14) to attain first order, that is, we set $\omega_2 = 1-\omega_1$. With these choices, we get a family of ASIRK-2A methods with 3 parameters,

$$\mathscr{B} = \begin{pmatrix} 0 & 0 \\ \omega_1+\gamma_1 & 0 \end{pmatrix}, \quad \mathscr{C} = \begin{pmatrix} \lambda_1 & 0 \\ \omega_1 & 1-\omega_1 \end{pmatrix}, \quad \omega = \begin{pmatrix} \omega_1 \\ 1-\omega_1 \end{pmatrix}. \quad (44)$$

To obtain second order, we impose conditions (15) and (16) to methods (44) obtaining

$$(1-\omega_1)(\omega_1+\gamma_1)=\tfrac{1}{2},$$
$$\omega_1\lambda_1+1-\omega_1=\tfrac{1}{2},$$

whose solution is given by

$$\lambda_1=\frac{2\omega_1-1}{2\omega_1},\quad \gamma_1=\frac{-2\omega_1^2+2\omega_1-1}{2(\omega_1-1)}. \quad (45)$$

In this way, we get a $\omega_1$-familiy of second order methods. It is not possible to achieve third order. However, the additional conditions for uniform convergence (29) are already satisfied because the last row of $\mathscr{C}$ is equal to $\omega$.

Now, from the stiff accuracy analysis (31), neither the term $\varepsilon h$ in the $u$ variable nor the term $\varepsilon^2/h$ in the $v$ variable can be cancelled. Consequently, the difference between the exact solution of problem (30) and the numerical solution given by methods (44)-(45) is of the form

$$u_1(h,\varepsilon)-u(h,\varepsilon)=\mathcal{O}(h^3,h\varepsilon),\quad v_1(h,\varepsilon)-v(h,\varepsilon)=\mathcal{O}(h^3,h\varepsilon,\varepsilon^2/h).$$

In order to remove the negative powers in $h$ we need more parameters in the scheme and for this reason we consider $s=3$.

*Case s=3.* Like in the previous case, we set $\lambda_3=\omega_3$ and we impose the first order condition (14), that is, $\omega_3=1-\omega_1-\omega_2$. With these choices, we get a family of ASIRK-3A methods with 6 parameters,

$$\mathscr{B}=\begin{pmatrix}0 & 0 & 0\\ \omega_1+\gamma_1 & 0 & 0\\ \omega_1 & \omega_2+\gamma_2 & 0\end{pmatrix},\quad \mathscr{C}=\begin{pmatrix}\lambda_1 & 0 & 0\\ \omega_1 & \lambda_2 & 0\\ \omega_1 & \omega_2 & 1-\omega_1-\omega_2\end{pmatrix},\quad \omega=\begin{pmatrix}\omega_1\\ \omega_2\\ 1-\omega_1-\omega_2\end{pmatrix}. \quad (46)$$

Order 2 is assured if equations (15)-(16) are satisfied. For methods (46) these equations are

$$\omega_2(\omega_1+\gamma_1)+(-\omega_1-\omega_2+1)(\omega_1+\omega_2+\gamma_2)=\frac{1}{2}, \quad (47a)$$

$$\omega_1\omega_2+\omega_1\lambda_1-\omega_1+\omega_2\lambda_2-\omega_2+1=\frac{1}{2}. \quad (47b)$$

It can be checked that, with the 6 parameters in (46), it is not possible to achieve order three.

For second order ASIRK-3A methods of the form (46), the additional order conditions (29) are satisfied; observe that in (46), the last row of $\mathscr{C}$ is equal to $\omega$.

Furthermore, in (46) we set $\lambda=\lambda_1=\lambda_2$. In this way, if function $g$ is linear, or we use Newton-like methods for solving the nonlinear systems, the same LU-factorization can be used when the first two internal derivatives are computed (see equation (43)).

From the stiff accuracy analysis (31), cancellation of the term $\varepsilon h$ in the $u$ variable, and the term $\varepsilon^2/h$ in the $v$ variable, leads to the equations

$$\frac{\omega_2\lambda(\omega_1+\gamma_1)+(\omega_1+\omega_2-1)(\omega_1(\omega_2+\gamma_2)-\lambda(\omega_1+\omega_2+\gamma_2))}{\lambda^2}=1, \quad (48)$$

$$\frac{(\omega_1-\lambda)(\omega_2-\lambda)}{\lambda^2(\omega_1+\omega_2-1)}=0. \quad (49)$$

If (48) and (49) hold, the difference between the numerical and the exact solution of problem (30) is of the form

$$u_1(\varepsilon,h)-u(\varepsilon,h)=\mathcal{O}(h^3,h^2\varepsilon),\quad v_1(\varepsilon,h)-v(\varepsilon,h)=\mathcal{O}(h^3,h\varepsilon). \quad (50)$$

In order to impose condition (49), we take $\lambda = \omega_1$, and substitute this value in the other three equations, namely, (47a), (47b) and (48). Thus, we get the system

$$\omega_2(\omega_1 + \gamma_1) + (-\omega_1 - \omega_2 + 1)(\omega_1 + \omega_2 + \gamma_2) = \frac{1}{2},$$

$$\omega_2 = \frac{-2\omega_1^2 + 2\omega_1 - 1}{2(2\omega_1 - 1)},$$

$$\omega_1^2 = \omega_2 \gamma_1,$$

whose solution is

$$\gamma_2 = \frac{4\omega_1^4 + 4\omega_1^3 - 12\omega_1^2 + 6\omega_1 - 1}{2(4\omega_1^3 - 10\omega_1^2 + 6\omega_1 - 1)},$$

$$\omega_2 = \frac{-2\omega_1^2 + 2\omega_1 - 1}{2(2\omega_1 - 1)}, \tag{51}$$

$$\gamma_1 = -\frac{2(2\omega_1^3 - \omega_1^2)}{2\omega_1^2 - 2\omega_1 + 1}.$$

In this way, we get a $\omega_1$-family of second order methods such that the implicit method is *L*-stable, provided that it is *A*-stable. Furthermore, the methods satisfy the additional conditions (29), they can be implemented by using just three memory registers and, for the model problem (30), the errors are given by (50).

The next step is to choose the parameter $\omega_1$ in (51). We will follow two different approaches to get it. Firstly we minimize the local truncation error, obtaining the scheme named *ASIRK-LSe(3,2)*. Secondly we optimize the stability region (26), obtaining the scheme named *ASIRK-LSs(3,2)*.

**ASIRK-LSe(3,2) method**

We have denoted by ASIRK-LSe(3,2) the method in (51) that minimizes the sum of the absolute values of the coefficients in the leading truncation error term. The optimum value is 0.218 and it is obtained for $\omega_1 = 0.1547777646705979$. By using the rational approximation $\omega_1 = 3/20$, we obtain the ASIRK-LSe(3,2) method

$$\mathscr{B} = \begin{pmatrix} 0 & 0 & 0 \\ \frac{573}{2980} & 0 & 0 \\ \frac{3}{20} & \frac{98}{89} & 0 \end{pmatrix}, \quad \mathscr{C} = \begin{pmatrix} \frac{3}{20} & 0 & 0 \\ \frac{3}{20} & \frac{3}{20} & 0 \\ \frac{3}{20} & \frac{149}{280} & \frac{89}{280} \end{pmatrix}, \quad \omega = \begin{pmatrix} \frac{3}{20} \\ \frac{149}{280} \\ \frac{89}{280} \end{pmatrix}. \tag{52}$$

The stability function (25) for this methods is a rational function that can be easily obtained from this formula

$$R(z_1, z_2) = \frac{\det(I_3 - z_1\mathscr{B} - z_2\mathscr{C} + (z_1 + z_2)e\omega^t)}{\det(I_3 - z_1\mathscr{B} - z_2\mathscr{C})}. \tag{53}$$

For method (52), it is

$$R(z_1, z_2) = \frac{59600(107z_2 + 280) + 59600(107z_2 + 280)z_1 + (1003731z_2 + 8344000)z_1^2 + 1123080z_1^3}{149(280 - 89z_2)(20 - 3z_2)^2}.$$

**ASIRK-LSs(3,2) method**

In this case we choose the parameter $\omega_1$ in the family (51) so as to maximize the area of the stability region $S_1$ in (26). In figure 1 we show the upper half of the stability region for different values of the parameter $\omega_1$, including the one that gives method (52). After analyzing these regions, we choose $\omega_1 = 0.14$ to get the ASIRK-LSs(3,2) scheme

$$\mathscr{B} = \begin{pmatrix} 0 & 0 & 0 \\ \frac{8407}{47450} & 0 & 0 \\ \frac{7}{50} & \frac{648}{599} & 0 \end{pmatrix}, \quad \mathscr{C} = \begin{pmatrix} \frac{7}{50} & 0 & 0 \\ \frac{7}{50} & \frac{7}{50} & 0 \\ \frac{7}{50} & \frac{949}{1800} & \frac{599}{1800} \end{pmatrix}, \quad \omega = \begin{pmatrix} \frac{7}{50} \\ \frac{149}{280} \\ \frac{599}{1800} \end{pmatrix}. \quad (54)$$

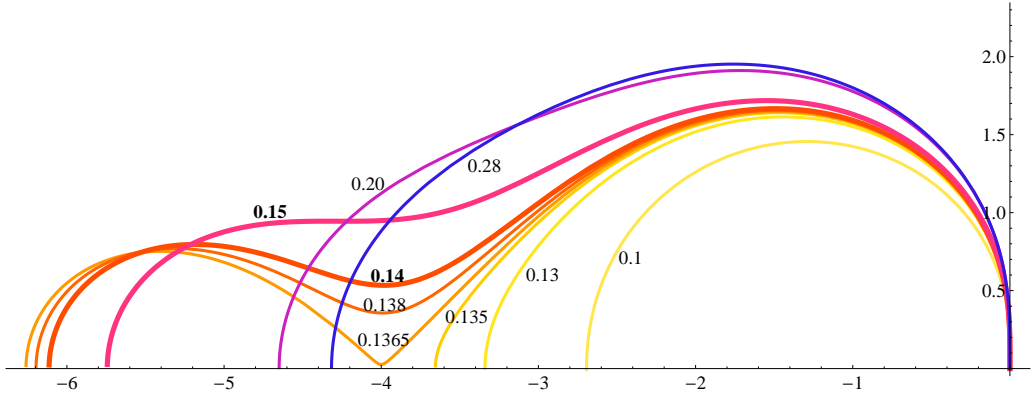

Figure 1: Upper half of the stability region $S_1$ for method (51) for different values of the parameter $\omega_1$: 0.1, 0.13, 0.135, 0.1365. 0.138. 0.14, 0.15, 0.20, 0.28. In particular, $\omega_1 = 0.15$ and $\omega_1 = 0.14$ correspond to methods (52) and (54), respectively.

For this method, the sum of the absolute values of the coefficients in the leading truncation error term is equal to 0.220; observe that this figure is close to the minimum value obtained for scheme (52). The stability function (25) for this scheme can be easily obtained from formula (53).

**Proposition 5.** *If we take $\lambda = \omega_2$ in (49), and we proceed in the same way, we get another family of second order methods whose coefficients are expressed in a more complicated way. The numerical experiments done, not shown in the paper, give slightly worse results than the ones obtained with the above family.*

## 5 Numerical Experiments

In this section we study the performance of the numerical schemes constructed in this paper and we compare them with other methods from the literature. More precisely, we consider the ASIRK-LSe(3,2) and ASIRK-LSs(3,2) methods given by (52) and (54), respectively, together with the ones we show below. In table 1 we summarize the properties of all them.

*ASIRK-LS(3,2):* In order to check the relevance of conditions (48)-(49), we have considered a second order method within the family (46) that does not satisfy them. The coefficients of this method are

$$\mathscr{B} = \begin{pmatrix} 0 & 0 & 0 \\ 0.679529 & 0 & 0 \\ 0.429529 & 0.591085 & 0 \end{pmatrix}, \quad \mathscr{C} = \begin{pmatrix} 0.1 & 0 & 0 \\ 0.429529 & 0.1 & 0 \\ 0.429529 & 0.241085 & 0.329385 \end{pmatrix}, \quad \omega = \mathscr{C}^t e_3. \quad (55)$$

As conditions (48)-(49) are not satisfied, the difference between the numerical and the exact solution of problem (30) is

$$u_1(\varepsilon, h) - u(\varepsilon, h) = \mathscr{O}(h^3, h\varepsilon), \qquad v_1(\varepsilon, h) - v(\varepsilon, h) = \mathscr{O}(h^3, h\varepsilon, \varepsilon^2/h).$$

| Method | Order | Additional cond. ((28) or (29)) | Stiff accuracy $(u,v)$ problem (30) | Registers of memory |
|---|---|---|---|---|
| ASIRK-LSe(3,2) (52) | 2 | Yes | $(\mathcal{O}(h^3,h^2\varepsilon), \mathcal{O}(h^3,h\varepsilon))$ | 3 |
| ASIRK-LSs(3,2) (54) | 2 | Yes | $(\mathcal{O}(h^3,h^2\varepsilon), \mathcal{O}(h^3,h\varepsilon))$ | 3 |
| ASIRK-LS(3,2) (55) | 2 | Yes | $(\mathcal{O}(h^3,h\varepsilon), \mathcal{O}(h^3,h\varepsilon,\varepsilon^2/h))$ | 3 |
| Zhong's method (56) | 2 (3) | No | $(\mathcal{O}(h^3,h\varepsilon), \mathcal{O}(h^2,h\varepsilon,\varepsilon^2/h))$ | $\geq 4$ |
| IMEX-SSP2(3,3,2) (57) | 2 | No | $(\mathcal{O}(h^3,h^2\varepsilon), \mathcal{O}(h^2,h\varepsilon))$ | $\geq 4$ |

Table 1: Properties of methods used in the numerical experiments: convergence order, stiff accuracy and number of memory registers required

This scheme satisfies the additional conditions (28). Furthermore, the implicit part is *L*-stable, the last row of $\mathscr{C}$ is equal to $\omega$, and it can be implemented in three memory registers.

*Zhong's method:* We consider the ASIRK-3A method [28]

$$\begin{aligned}
&b_{21} = \tfrac{8}{7}, \quad b_{31} = \tfrac{71}{252}, \quad b_{32} = \tfrac{7}{36}, \\
&c_{21} = 0.306727, \; c_{31} = 0.45, \quad c_{32} = -0.263111, \\
&c_{11} = 0.485561, \; c_{22} = 0.951130, \; c_{33} = 0.189208, \\
&\omega_1 = \tfrac{1}{8}, \quad \omega_2 = \tfrac{1}{8}, \quad \omega_3 = \tfrac{3}{4}.
\end{aligned} \quad (56)$$

This method satisfies the order equations (14)-(20) and (23). Consequently it is a third order method for problems such that the matrices $f_y$ and $g_y$ commute; however, for a general problem, it is a method of order 2. The implicit part in the scheme (56) is *L*-stable. However, $\omega^t \mathscr{C}^{-1}(\mathscr{C}e)^2 \neq 1$ and thus method (56) does not satisfy all the additional conditions (29). This method does not satisfy equations (48)-(49) either and, consequently, the difference between the numerical and the exact solution of problem (30) is

$$u_1(\varepsilon,h) - u(\varepsilon,h) = \mathcal{O}(h^3, h\varepsilon), \qquad v_1(\varepsilon,h) - v(\varepsilon,h) = \mathcal{O}(h^2, h\varepsilon, \varepsilon^2/h).$$

As there is not any pattern on the coefficients of the method, at least four memory registers are needed to implement it.

*IMEX-SSP2(3,3,2):* We consider the IMEX Runge-Kutta method [23]

$$\begin{array}{c|ccc}
0 & 0 & 0 & 0 \\
\tfrac{1}{2} & \tfrac{1}{2} & 0 & 0 \\
1 & \tfrac{1}{2} & \tfrac{1}{2} & 0 \\
\hline
& \tfrac{1}{3} & \tfrac{1}{3} & \tfrac{1}{3}
\end{array}
\qquad
\begin{array}{c|ccc}
\tfrac{1}{4} & \tfrac{1}{4} & 0 & 0 \\
\tfrac{1}{4} & 0 & \tfrac{1}{4} & 0 \\
1 & \tfrac{1}{3} & \tfrac{1}{3} & \tfrac{1}{3} \\
\hline
& \tfrac{1}{3} & \tfrac{1}{3} & \tfrac{1}{3}
\end{array}
\qquad (57)$$

This scheme is constructed in the context of strong-stability-preserving (SSP) schemes [9, 11, 18, 21, 24, 25] and references therein. It is a second order method with an *L*-stable implicit part. The coefficients of this scheme do not satisfy the additional conditions (28) for stiff accuracy since $b^t A^{-1} \tilde{A} \tilde{c} = 1/2$. Furthermore, the difference between the numerical and the exact solution of problem (30) is

$$u_1(\varepsilon,h) - u(\varepsilon,h) = \mathcal{O}(h^3, h^2\varepsilon), \qquad v_1(\varepsilon,h) - v(\varepsilon,h) = \mathcal{O}(h^2, h\varepsilon). \qquad (58)$$

Like Zhong's method (56), implementation of method (57) requires at least four memory registers.

## 5.1 Methods with a stiffness parameter $\varepsilon$

For our numerical tests, we have considered three problems containing a stiffness parameter $\varepsilon$: a simple prototype of additive system [22], the van der Pol equation [2] and the Broadwell problem [5]. For all of them, the matrices $f_y$ and $g_y$ do not commute and, consequently, all the methods are second order schemes. We will test the convergence rates when $\varepsilon$ varies and the efficiency of the methods for a given value of $\varepsilon$.

We have integrated these problems with different initial data: inconsistent initial values (IC$_{\text{InVal}}$), consistent initial values (C$_{\text{InVal}}$), and consistent *well prepared* initial values (WP$_{\text{InVal}}$). In the first case the problem presents an initial layer. On the contrary, *well prepared* initial values ensure a smooth solution for any value of $\varepsilon$.

*Problem 1:* [22] We consider the simple prototype of stiff system of the form

$$\begin{cases} u' = -v, \\ v' = u + \dfrac{1}{\varepsilon}(e(u) - v), \end{cases} \tag{59}$$

with $e(u) = \sin u$. We can choose arbitrarily the initial value $u(0)$, but there is no freedom in the choice of $v(0)$. If we consider inconsistent initial values, then the solution presents an initial layer in the $v$ component when $\varepsilon \to 0$. We have chosen $u(0) = \pi/2$ and we have integrated the problem in $[0, 1]$, with time step $h = 0.05$, assuming different initial data $v(0)$. Consistent initial values are obtained if we take $v(0) = \sin(u(0)) = 1$. By adding a small perturbation $\delta$, we get inconsistent initial values; in the numerical experiments we have taken $\delta = 0.05$. Finally we have also considered *well prepared* initial data to obtain a smooth solution.

$$\begin{align}
\text{C}_{\text{InVal}}: &\qquad v(0) = 1, \\
\text{IC}_{\text{InVal}}: &\qquad v(0) = 1 + \delta, \\
\text{WP}_{\text{InVal}}: &\qquad v(0) = 1 + \frac{\pi}{2}\varepsilon - \frac{\pi}{2}\varepsilon^3.
\end{align} \tag{60}$$

*Problem 2:* Another example of stiff problem of the form (2) is the van der Pol equation

$$\begin{cases} y' = z, \\ \varepsilon z' = (1 - y^2)z - y. \end{cases} \tag{61}$$

This model has also been considered in [3] in the context of uniform convergence. For this problem we can choose the initial value $y(0)$ arbitrarily, but there is no freedom in the choice of $z(0)$. So, we have taken $y(0) = 2$, and we have integrated the problem in $[0, 0.55139]$, with $h = 0.05$, assuming different initial data in a similar way as in the previous problem.

$$\begin{align}
\text{C}_{\text{InVal}}: &\qquad z(0) = -2/3, \\
\text{IC}_{\text{InVal}}: &\qquad z(0) = -2/3 + \delta, \\
\text{WP}_{\text{InVal}}: &\qquad z(0) = -\frac{2}{3} + \frac{10}{81}\varepsilon - \frac{292}{2187}\varepsilon^2 - \frac{1814}{19683}\varepsilon^3.
\end{align} \tag{62}$$

*Problem 3:* We consider the Broadwell model [5], a hyperbolic system with relaxation. We have integrated this problem in the fluid dynamical moment variables $\rho$, $m$ and $u$. In these variables, the Broadwell equations can be written as

$$\begin{align}
\partial_t \rho + \partial_x m &= 0, \\
\partial_t m + \partial_x z &= 0, \\
\partial_t z + \partial_x m &= \frac{1}{2\varepsilon}\left(\rho^2 + m^2 - 2\rho z\right).
\end{align} \tag{63}$$

A conservative spatial discretization applied to (63) gives (see [5, p. 251] for details)

$$\begin{aligned}
\partial_t \rho_j + \frac{m_{j+1} - m_{j-1}}{2\Delta x} - \frac{z_{j+1} - 2z_j + z_{j-1}}{2\Delta x} &= 0, \\
\partial_t m_j + \frac{z_{j+1} - z_{j-1}}{2\Delta x} - \frac{m_{j+1} - 2m_j + m_{j-1}}{2\Delta x} &= 0, \\
\partial_t z_j + \frac{m_{j+1} - m_{j-1}}{2\Delta x} - \frac{z_{j+1} - 2z_j + z_{j-1}}{2\Delta x} &= \frac{1}{2\varepsilon} \left( \rho_j^2 + m_j^2 - 2\rho_j z_j \right).
\end{aligned} \quad (64)$$

We have considered periodic boundary conditions over the domain $[-1,1]$ and we have integrated the problem in the time interval $[0, 0.5]$, with $\Delta x = 0.2$ and $\Delta t = 0.05$. We have taken the following initial values: $\rho(x,0) = \rho_0(x)$ and $m(x,0) = m_0(x)$, with

$$\rho_0(x) = 1 + a_p \sin\left(\tfrac{2\pi x}{L}\right), \qquad m_0(x) = \rho_0(x)\left(\tfrac{1}{2} + a_u \sin\left(\tfrac{2\pi x}{L}\right)\right),$$

where $a_p = 0.3$, $a_u = 0.1$, and $L = 2$ is the length of the domain. There is no freedom in the choice of $z(x,0) = z_0(x)$, and we have considered the following initial data

$$\begin{aligned}
\text{C}_{\text{InVal}}: & \qquad z_0(x) = z_E(\rho_0(x), m_0(x)), \\
\text{IC}_{\text{InVal}}: & \qquad z_0(x) = z_E(\rho_0(x), m_0(x)) + \delta, \\
\text{WP}_{\text{InVal}}: & \qquad z_0(x) = z_E(\rho_0(x), m_0(x)) + \varepsilon z_1(\rho_0(x), m_0(x)),
\end{aligned} \quad (65)$$

where

$$z_E(\rho, m) = \frac{\rho^2 + m^2}{2\rho}, \qquad z_1(\rho, m) = \frac{H(\rho, m)}{2\rho},$$

and

$$H(\rho, m) = \left(-1 + \partial_\rho z_E + (\partial_m z_E)^2\right) \partial_x m_0 + \partial_\rho z_E \, \partial_m z_E \, \partial_x \rho_0.$$

#### Covergence rate for different values of $\varepsilon$

For each method and problem, we compute the convergence rates for a wide range of values of the parameter $\varepsilon$; in our tests we consider $\varepsilon = 10^{-j}$, $j = 0, 1, 2, \ldots, 6$. We try to check whether the convergence is uniform in $\varepsilon$, particularly in the intermediate regime. In order to show the convergence rates, for each value of $\varepsilon$ we compute $E_h$ and $E_{h/2}$, an estimation of the relative $L_\infty$-global error with stepsizes $h$ and $h/2$, respectively. With these values we compute the convergence rate in the standard way.

In figures 2, 3 and 4 we show the convergence rates versus the stiff parameter $\varepsilon$ for problems (59), (61) and (64), respectively, when methods ASIRK-LSe(3,2), ASIRK-LSs(3,2), ASIRK-LS(3,2), Zhong's method and IMEX-SSP2(3,3,2) method, given by equations (52), (54), (55), (56) and (57), respectively, are used.

The best numerical results are obtained for the two new methods constructed in this paper named ASIRK-LSe(3,2) and ASIRK-LSs(3,2). For WP$_{\text{InVal}}$ initial conditions, they show second order uniform convergence for the Broadwell problem (see figure 4) and almost second order of convergence for the rest of the problems. For problem (59), in the middle stiff regime, that is, when $h$ is close to $\varepsilon$, the orders for methods (52) and (54) are decreased to 1.68 and 1.71, respectively; for the van der Pol problem (see figure 3), the orders are reduced to 1.80 and 1.82, respectively. When C$_{\text{InVal}}$ initial values are used, the results are quite similar, although in this case order reduction is more evident in some cases. Methods (52) and (54) show second order uniform convergence for the Broadwell problem but, for the van der Pol equation, in the middle stiff regime, the order is reduced to 0.52 and 0.53, respectively; the reduction is larger for problem (59), where we obtain values of $-1.12$ and $-1.09$, respectively. For IC$_{\text{InVal}}$ initial conditions, for most of the problems, the best numerical results are obtained with the new robust schemes constructed in this paper.

Method ASIRK-LS(3,2) in (55) does not fulfill the stiff accuracy conditions (48) and (49). When WP$_{\text{InVal}}$ initial values are used, this method shows uniform second order for Broawell model but not for the van der Pol equation where, for some values of $\varepsilon$, the order is dramatically decreased. For C$_{\text{InVal}}$ initial conditions, the order decreases up to 1.35 in the middle stiff regime when Broadwell problem is integrated. However, for the rest of the problems, the error is dramatically dropped. The numerical results are worst when IC$_{\text{InVal}}$ initial conditions are used.

Zhong's method (56) does not fulfill neither the stiff accuracy conditions (48) and (49) nor the additional conditions (29). Even for WP$_{\text{InVal}}$ initial values, a poor result is obtained for the van der Pol problem; observe that this method applied to the Broadwell model shows third order rate in some components when small values of $\varepsilon$ are considered. The numerical results are worst when IC$_{\text{InVal}}$ initial conditions are used.

IMEX-SSP2(3,3,2) scheme (57) does not fulfill the additional condition (28) and, for problem (30), the errors do not contain negative powers of the step size $h$ (see equation (58)). For WP$_{\text{InVal}}$ initial conditions, we observe uniform second order of convergence for problem (61) whereas, for the van der Pol problem and the Broadwell problem the order is reduced to 1.68 and 1.58, respectively. Taking into account all the numerical experiments done, the results with IMEX-SSP2(3,3,2) method are worse than the ones obtained with the two new schemes constructed in this paper.

With regard to the behaviour when $\varepsilon$ tends to zero, the new ASIRK methods (52) and (54) attain second order for all the problems and all the initial conditions considered. The same order is obtained for IMEX-SSP2(3,3,2) scheme (57). For these three schemes, the errors for problem (30) do not contain negative powers of the step size $h$. On the other hand, for methods like ASIRK-LS(3,2) method (55) and Zhong's method (56), whose errors contain negative powers of $h$, first order is observed when IC$_{\text{InVal}}$ initial values are used for problem (59) and for the van der Pol problem.

### *Efficiency*

Next, for a given value of the parameter $\varepsilon$, we have computed the errors when different step sizes are used. In figures 5 and 6, we show the number of steps versus the error when methods ASIRK-LSe(3,2), ASIRK-LSs(3,2), ASIRK-LS(3,2), Zhong's method and IMEX-SSP2(3,3,2) method are used for problems (59) and (61), respectively. We have taken $\varepsilon = 1.e-5$ and $\varepsilon = 1.e-3$ for problem (59) and (61), respectively; in both cases, we consider C$_{\text{InVal}}$ initial values.

In problem (59), for both variables the smallest errors are obtained with methods ASIRK-LSe(3,2) and ASIRK-LSs(3,2). Furthermore, the results for variable $v$ (see figure 5, right) show that the new methods constructed in this paper performs much better than the ASIRK method by Zhong or the IMEX-SSP2(3,3,2) scheme. Among the new methods constructed, the schemes ASIRK-LSe(3,2) and ASIRK-LSs(3,2) behave in a similar way, whereas method ASIRK-LS(3,2) shows larger errors.

Similar results are obtained in problem (61). In this case, for variable $y$, the errors with the new ASIRK methods in this paper are slighter smaller but, for variable $z$, they are smaller than the ones obtained with the ASIRK method by Zhong or the IMEX-SSP2(3,3,2) scheme.

## 5.2 Population model

Next, in order to show the performance of the methods in this paper for problems without a stiffness parameter, we consider the problem in [15] that models the evolution of a population density $P$,

$$\frac{\partial P(t,x)}{\partial t} = f(t,x) + b(x,P(t,x)) - r_d P(t,x) + d\frac{\partial^2 P(t,x)}{\partial x^2}, \tag{66}$$

where $x \in [0,1]$ and $t \geq 0$. We assume periodic boundary conditions and zero initial condition, $P(0,t) = 0$, $x \in [0,1]$. The birth rate $b(x,P)$ is given by

$$b(x,P) = r_b(x)\frac{\varepsilon}{\varepsilon + P}, \qquad \varepsilon = 0.005,$$

where

$$r_b(x) = \begin{cases} 1, & \text{if } x \in [0,1/2], \\ 100, & \text{if } x \in (1/2,1]. \end{cases}$$

The death rate is set to one, i.e., $r_d = 1$, and the forcing term will be described below. In order to discretize (66), a spatial grid $\Delta x = 1/100$ is taken, and second order centered differences are used to discretize the diffusion term. The function $f(t,x)$ is taken to be zero for all times $t \neq 0$, and for each grid point $x_i = i\Delta x$, $f(0,x_i)$ is a random value in the interval $[0.8, 1.2]$. In the time stepping process, the diffusion term is integrated implicitly and the rest of the terms are integrated explicitly. In the numerical experiments we have considered $d = 0.02$, $d = 0.04$.

In figure 7 we show the numerical solution at $t_{end} = 10$ obtained with ASIRK-LSe(3,2), ASIRK-LS(3,2) and Zhong's methods with step size $h = 20/9$. For comparisons, a reference solution has been computed with Matlab with a small tolerance. We see that the best numerical solution is given by ASIRK-LS(3,2) scheme, followed by ASIRK-LSe(3,2) method; the worst results are obtained with Zhong's method. The numerical results obtained with ASIRK-LSs(3,2) scheme (not shown in this paper), are similar to the ones obtained with ASIRK-LSe(3,2) method. Interestingly, methods with the best performance, namely ASIRK-LSe(3,2), ASIRK-LSs(3,2) and ASIRK-LS(3,2), satisfy all the additional conditions (29).

# 6 Concluding remarks

In this paper we have constructed two robust second order 3-stage ASIRK methods, whose coefficients are given in (52) and (54), that can be implemented with just three memory registers. The new methods incorporate accuracy, stiff accuracy and linear stability properties. The numerical experiments with different problems and different initial conditions show that they outperform some other methods of the literature that do not share these properties.

The numerical results for the two new ASIRK schemes are quite similar although, for most of the problems and initial values, method (54) gives slightly better results. This behavior can be explained because, the stability region for method (54) is slightly larger than the one for method (52) (see figure 1), but the size of the leading truncation error terms is quite similar for both schemes: 0.218 for method (52), and 0.220 for method (54).

In order to improve the numerical results when $C_{InVal}$ are used, some additional conditions could be imposed to the scheme if we would have more free parameters to construct it. These additional parameters can be obtained by increasing the number of stages of the method. Observe that, for the class of low-storage methods considered in this paper, this action does not increase the number of memory registers required to implement the scheme.

Some important problems, e.g., the resolution of nonlinear systems with low-storage strategies, remain open. They will be a topic of future research.

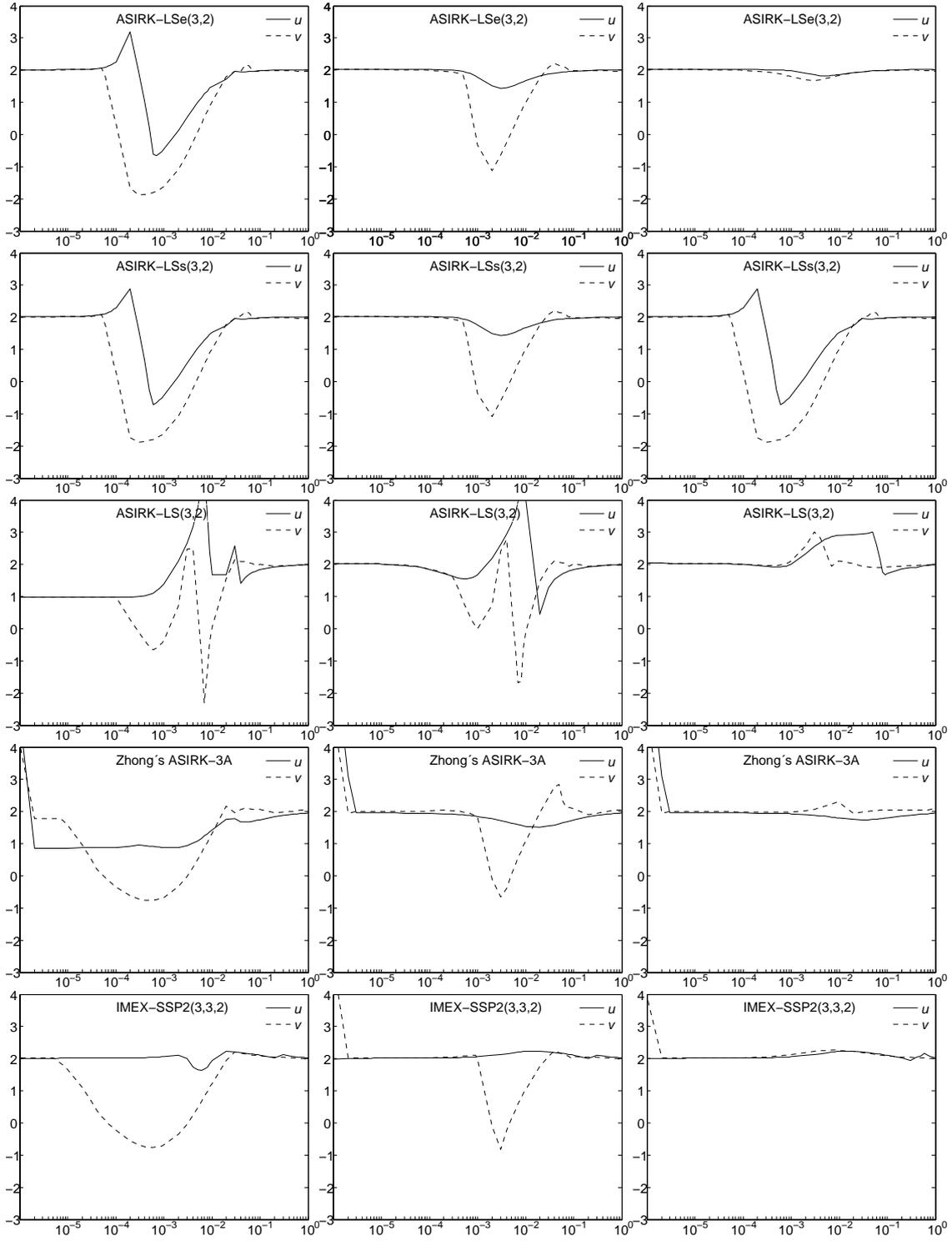

Figure 2: Convergence rates of $u$ (continuous line) and $v$ (dashed line) vs. stiffness parameter $\varepsilon$ for problem (59) for stepsize $h = 0.05$.
From top to bottom: ASIRK-LSe(3,2) (52), ASIRK-LSs(3,2) (54), ASIRK-LS(3,2) (55), Zhong's method (56) and IMEX-SSP2(3,3,2) (57). From left to right: $IC_{InVal}$, $C_{InVal}$ and $WP_{InVal}$ initial values (60)

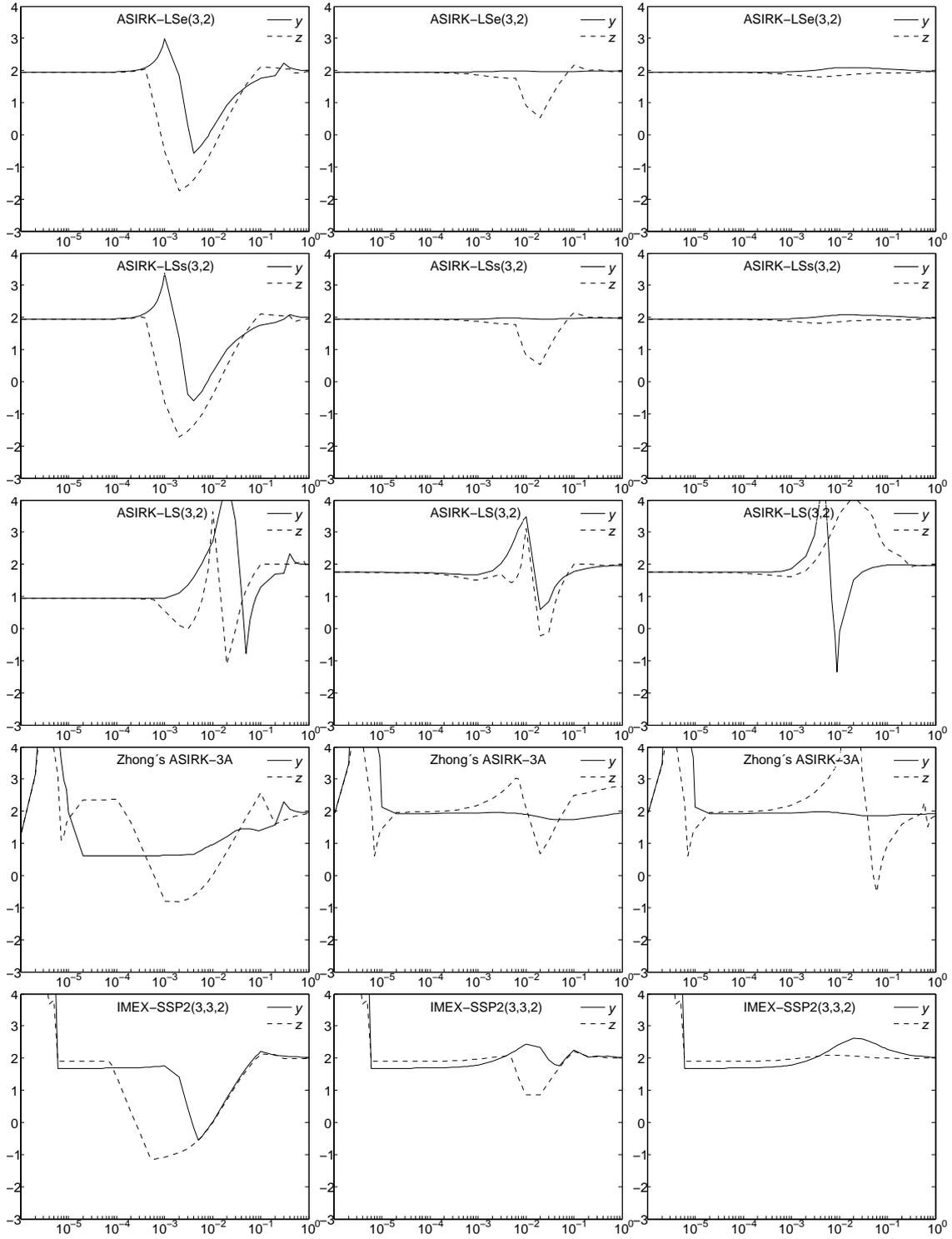

Figure 3: Convergence rates of $y$ (continuous line) and $z$ (dashed line) vs. stiffness parameter $\varepsilon$ for van der Pol problem (61) for stepsize $h = 0.05$.
From top to bottom: ASIRK-LSe(3,2) scheme (52), ASIRK-LSs(3,2) scheme (54), ASIRK-LS(3,2) scheme (55), Zhong's method (56) and IMEX-SSP2(3,3,2) (57). From left to right: $IC_{InVal}$, $C_{InVal}$ and $WP_{InVal}$ initial values (62)

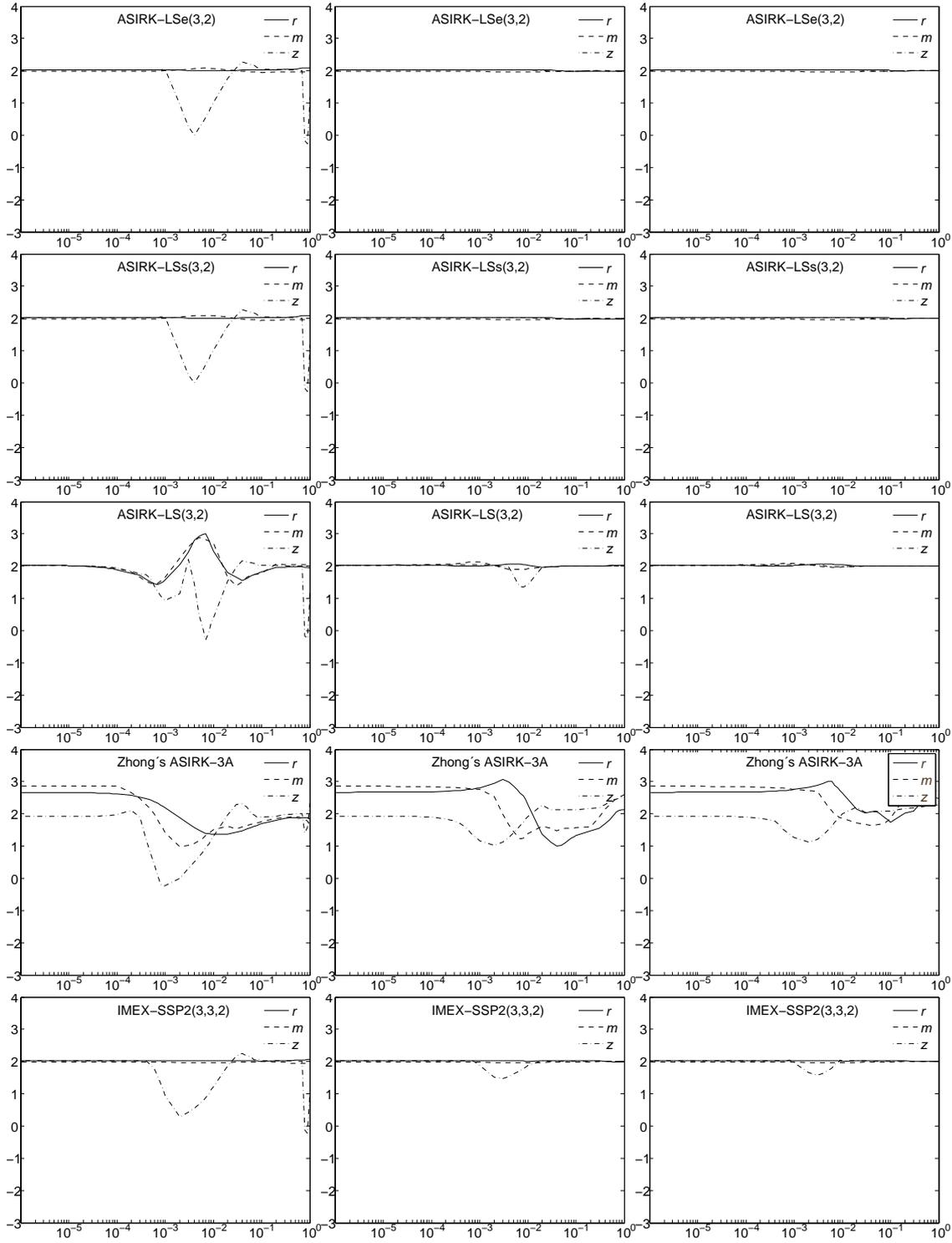

Figure 4: Convergence rates of $r$ (continuous line), $m$ (dashed line) and $z$ (dotted-dashed line) vs. stiffness parameter $\varepsilon$ for Broadwell problem (64) for stepsize $h = 0.05$.
From top to bottom: ASIRK-LSe(3,2) scheme (52), ASIRK-LSs(3,2) scheme (54), ASIRK-LS(3,2) scheme (55), Zhong's method (56) and IMEX-SSP2(3,3,2) (57). From left to right: IC$_{\text{InVal}}$, C$_{\text{InVal}}$ and WP$_{\text{InVal}}$ initial values (65)

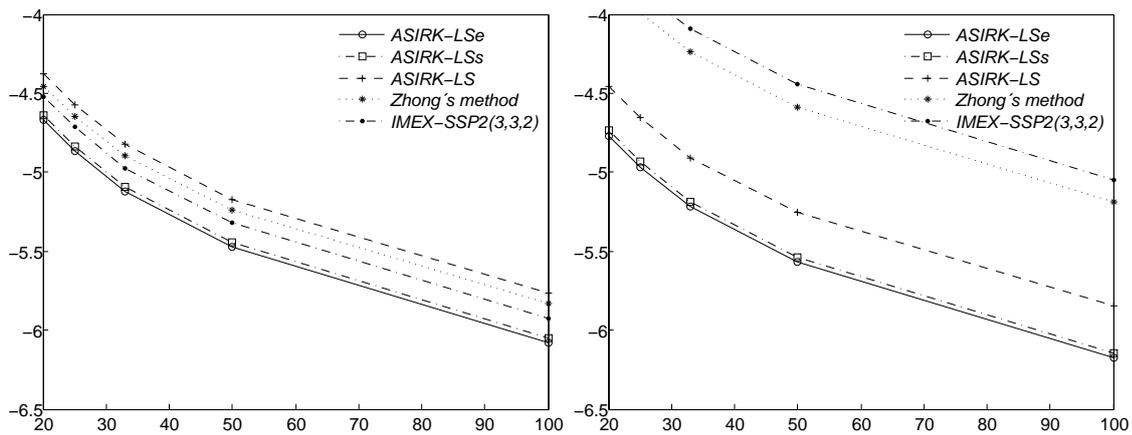

Figure 5: Error vs number of steps in variables $u$ (left) and $v$ (right) for problem (59) with $C_{\text{InVal}}$ initial values. $\varepsilon = 1e-5$.

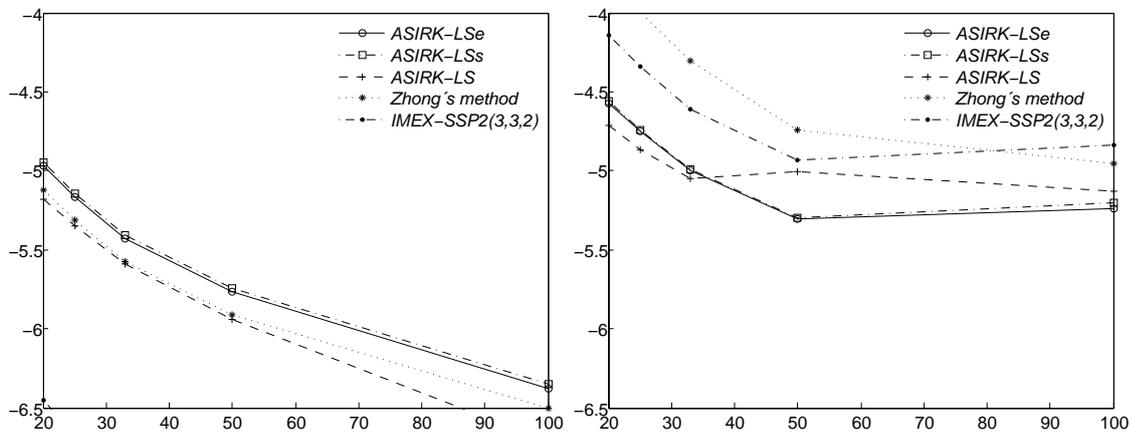

Figure 6: Error vs number of steps in variables $y$ (left) and $z$ (right) for problem (61) with $C_{\text{InVal}}$ initial values. $\varepsilon = 1e-3$.

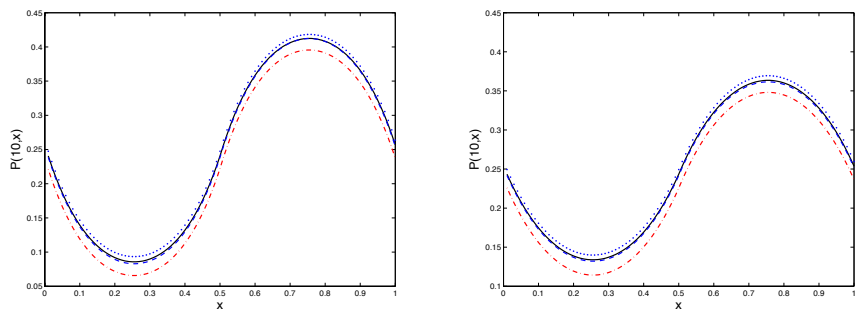

Figure 7: Problem (66). Numerical solution at $t_{end} = 10$ with step size $h = 20/9$ for $d = 0.02$ (left) and $d = 0.04$ (right): reference solution (solid line), ASIRK-LSe(3,2) (dotted line), ASIRK-LS(3,2) (dashed line) and Zhong's method (dash-dot line).